\documentclass[11pt]{article}
\usepackage{amsmath,amsfonts,amssymb,latexsym,amsbsy, bbm, theorem,enumerate,color}
\usepackage{graphicx, graphics}
\usepackage{float,color,fancybox,shapepar,setspace,hyperref}
\usepackage[affil-it]{authblk}
\usepackage{caption}
\usepackage{subfigure}
\usepackage{psfrag}
\usepackage{epstopdf}
\newtheorem{thm}{Theorem}[section]
\newtheorem{conj}[thm]{Conjecture}
\newtheorem{cor}[thm]{Corollary}
\newtheorem{lem}[thm]{Lemma}
\theorembodyfont{\rmfamily}
\def\pf{\medskip\noindent {\bf Proof.}~~}

\def\dfn#1{{\sl #1}}
\def\mytextindent#1{\indent\llap{#1\enspace}\ignorespaces}

\def\myitem{\par\hangindent\parindent\mytextindent}

\newcommand{\less}{\setminus}
\def\qed{ \hfill $\blacksquare$}

\begin{document}

\title{Antimagic orientations of graphs with given independence number}

\author{Zi-Xia Song$^{1,}$\thanks{Partially supported by the National Science Foundation under
grant   DMS-1854903.}\,, Donglei Yang$^{2,}$\thanks{This work was done  in part while the second author visited the University of Central Florida.}\,    and  Fangfang Zhang$^{3,}$\thanks{Corresponding author. This work was done while the third author visited the University of Central Florida as a visiting student.   The visit was 
supported by the Chinese Scholarship Council.  E-mail addresses:  Zixia.Song@ucf.edu (Z-X. Song); dlyang120@163.com (D. Yang);   and Fangfangzh@smail.nju.edu.cn (F. Zhang).}}
    
  \affil{ 
  { \small {$^1$Department  of Mathematics, University of Central Florida, Orlando, FL 32816, USA}} \\
  { \small {$^2$Department  of Mathematics, Shandong University, Jinan  250100, China}}\\ 
  { \small {$^3$Department of Mathematics, Nanjing University, Nanjing  210093, China}}\\ 
  }

 \date{ }
 
\maketitle

\begin{abstract}
Given a digraph $D$ with $m $ arcs and  a bijection $\tau: A(D)\rightarrow \{1, 2, \ldots, m\}$,   we say $(D, \tau)$ is  an
  \dfn{antimagic orientation}   of a graph $G$  if    $D$  is an orientation of $G$ and  no two vertices in $D$ have the same  vertex-sum under $\tau$,  where the  vertex-sum  of a vertex $u $  in $D$ under  $\tau$ is the sum of labels of all arcs entering $u$ minus the sum of labels of all arcs leaving $u$.   Hefetz, M\"{u}tze, and Schwartz in 2010 initiated the study of antimagic  orientations of graphs, and  conjectured that every connected graph  admits an antimagic orientation.   This conjecture seems  hard, and few related results are known. However, it has been verified
to be true for  regular graphs,   biregular bipartite graphs, and graphs with large maximum degree.
 In this paper,  we establish more evidence for the aforementioned conjecture by studying antimagic orientations of graphs $G$ with independence number at least $|V(G)|/2$ or at most four. We obtain several results. The method we develop in this paper may shed some light on attacking the aforementioned conjecture. 

\end{abstract}

\baselineskip 18pt
\section{Introduction}

In this paper,  all  graphs   are finite and  simple, and all multigraphs are finite and loopless.     For a graph $G$, we use $V(G)$, $E(G)$, $|G|$,  $e(G)$, $\Delta(G)$,  $\delta(G)$ and $\alpha(G)$ to denote the vertex set, edge set, number of vertices, number of   edges, maximum degree,   minimum degree, and independence number of $G$, respectively.  Given   sets $S\subseteq V(G)$ and $F\subseteq E(G)$,  we use        $G\less S$ to denote the subgraph    obtained from $G$ by deleting all vertices in $S$,   $G\less F$ the subgraph    obtained from $G$ by deleting all edges in $F$, and $G[S] $    the  subgraph    obtained from $G$ by deleting all vertices in $V(G)\less S$.   For  two disjoint sets $A, B\subseteq V(G)$,    $A$ is \dfn{complete} to $B$ in $G$  if each vertex in $A$ is adjacent to all vertices in  $B$, and \dfn{anti-complete} to $B$ in $G$  if no  vertex in $A$ is adjacent to any vertex  in  $B$.   
We simply say $a$ is \dfn{complete to} (resp. \dfn{anticomplete to}) $B$ when $A=\{a\}$.
For convenience, we use   $A \less B$ to denote  $A -B$; and $A \less b$ to denote  $A -\{b\}$ when $B=\{b\}$.  We use the convention   ``$S:=$'' to mean that $S$ is defined to be the right-hand side of the relation.  The degree
and neighborhood of a vertex $v$ in $G$ are denoted by $d_G(v)$ and $N_G(v)$, respectively. We
define  $N_G[v]:=N_G(v)\cup\{v\}$.   If no confusion arises, we omit the subscript $G$ in the above notation. 
 Let $H$ be a   connected graph.    The \dfn{distance}  $d_H(x, y)$ of two vertices  $x$ and $y$ is   the length of a shortest  $(x,y)$-path in $H$.  The \dfn{radius}  $r(H)$ of   $H$ is defined to be $\min_{x\in V(H)}\max_{y\in V(H)} d_H(x, y)$.  Given $X\subseteq V(H)$ and $v\in V(H)\less X$,  we define  $d_H(v, X):= min_{x\in X} d_H(v, x)$. 
  For a positive integer $k$, we write  $[k]$ for the set $\{1,2, \ldots, k\}$.  We use $P_n$,  $C_n$, $K_n$, $S_{n-1}$  and $W_{n-1}$ to denote the path,    cycle, complete graph, star and   wheel  on $n$ vertices, respectively. 

\medskip

  An \dfn{antimagic labeling} of a graph $G$ with $m$ edges is a bijection $\tau:  E(G) \rightarrow  [m]$ such that for any distinct vertices $u$ and $v$, the sum of labels on edges incident to $u$ differs from the sum of labels on edges incident  to $v$.   A graph   is \dfn{antimagic} if it has an antimagic labeling. Hartsfield and Ringel~\cite{NG} introduced antimagic labelings in 1990 and conjectured that every connected graph  other than $K_2$  is  antimagic.  The most recent  progress on this problem is a   result of   Eccles~\cite{E}, which states that   there exists an absolute constant $c_0$ such  that if   $G$  is a graph with average
degree at least $c_0$, and $G$ contains no isolated edge and at most one isolated vertex, then $G$ is
antimagic.  This    improves  a  result of Alon, Kaplan, Lev, Roditty, and Yuster~\cite{Alon}, which states that there exists an absolute   constant  $c$ such that  every   graph on $n$ vertices with minimum degree at least $c\log n$  is antimagic.    Hartsfield and Ringel's Conjecture has also been verified to be true for $d$-regular graphs with $d\ge2$ (see~\cite{KAM, CLPZ, DC, DYZ}),  and   graphs $G$  with $\Delta(G)\geq|G|-3\geq6$ by   Yilma \cite{Yilma}. For more information on antimagic labelings of graphs and related labeling problems, see the  recent  informative survey~\cite{JAG}.\medskip

Motivated by Hartsfield and Ringel's Conjecture, Hefetz, M\"{u}tze, and Schwartz~\cite{DTJ} introduced antimagic labeling  of digraphs. Let $D$ be a  digraph with $m$ arcs.  Let $\{a_1, \ldots, a_m\}$ be a set of  $m$ positive integers. For every bijection    $\tau: A(D)\rightarrow \{a_1,   \ldots, a_m\}$ and for each vertex  $u\in V(D)$, we define $s_{(D,\tau)}(u)$  to be
 the sum of labels of all arcs entering $u$ minus the sum of labels of all arcs leaving $u$ under $\tau$ when $u$ is not an isolated vertex in $D$, and $s_{(D,\tau)}(u):=0$ when $u$ is an isolated vertex in $D$;.   A bijection $\tau: A(D)\rightarrow [m]$ is     an \dfn{antimagic labeling} of $D$ if  $s_{(D,\tau)}(u)\neq s_{(D,\tau)}(v)$ for all   distinct vertices $u$ and $v$ in $D$.      A digraph $D$ is \dfn{antimagic} if it has an antimagic labeling.  
 We say $(D, \tau)$ is  an  \dfn{antimagic orientation}   of a graph $G$  if    $D$  is an orientation of $G$ and $\tau$ is       an antimagic labeling of $D$. 
 Hefetz, M\"{u}tze, and Schwartz~\cite{DTJ}  proved that   every orientation of       $S_n$ with $n\neq2$,  $K_n$ with $n\ne3$, and   $W_n$ with $n\ge3$ is antimagic, they further 
 asked whether it is true that    every orientation of any connected  graph, other than $K_3$ and $P_3$, is antimagic.   The same authors   proved an analogous result of Alon, Kaplan, Lev, Roditty, and Yuster~\cite{Alon}, which states that there exists an absolute constant $c$ such that every orientation of any  graph on $n$ vertices with minimum degree at least $c\log n$ is antimagic.  As pointed out in \cite{DTJ}, ``Proving that every orientation of such a graph is antimagic, however, seems rather difficult.''  As a relaxation of this problem, Hefetz, M\"{u}tze, and Schwartz~\cite{DTJ}  proposed the following conjecture.
 
 \begin{conj}[\cite{DTJ}]\label{conj}
Every connected graph admits an antimagic orientation.
\end{conj}
 Conjecture~\ref{conj} has been verified to be true for  odd regular graphs   \cite{DTJ},   disjoint union of cycles or  connected $2d$-regular graphs with $d\ge2$  by Li, Song, Wang, Yang, and Zhang   \cite{LSWY}, and   disconnected $2d$-regular graphs with $d\ge2$ by Yang  \cite{Y}.  
 A  bipartite graph $G$ with bipartition $(A, B)$  is \dfn{$(a,b)$-biregular} if   each vertex in $A$ has degree $a$ and each vertex in $B$ has degree $b$.   Shan and Yu~\cite{SY} recently proved that  every $(a, b)$-biregular bipartite graph admits an antimagic orientation.
      Very recently,  it has been proven in \cite{maxdegree} that Conjecture~\ref{conj} holds for connected graphs $G$ on $n\ge9$ vertices with $\Delta(G)\ge n-5$  or graphs with a dominating set of size two. 

\begin{thm}[\cite{maxdegree}]\label{dominating}
Let $G$ be a graph and let  $x,y\in V(G)$ be distinct  such that   $N[x]\cup N[y] =V(G)$ and $d(x)\ge d(y)$. If    $d(x)\ge4$    or $ N(x)\cap N[y] \ne\emptyset $, then $G$ admits an antimagic orientation.
\end{thm}

\begin{thm}[\cite{maxdegree}]\label{main} Let $G$ be a connected graph. 
\begin{enumerate}[\rm(i)] 
\item If  $\Delta(G)\geq|G|-3$, then $G$ has an antimagic orientation.
\item If  $\Delta(G)=|G|-t\ge4$, where   $t\in\{4,5\}$,  then   $G$ has an antimagic orientation.
\end{enumerate} 
\end{thm}
 \medskip
 
In this paper,   we establish more evidence for  Conjecture \ref{conj}   by studying antimagic orientations of  graphs $G$   with $\alpha(G)\ge |G|/2$ or $\alpha(G)\le 4$.   Theorem~\ref{pm} is a result from \cite{DTJ}. 

\begin{thm}[\cite{DTJ}]\label{pm} 
Let $G$ be a graph on $2n$ vertices that admits a perfect matching, and let $A$ be an independent set in $G$ with $|A|=n$. If $d(v)\ge3$ for every $v\in A$, then $G$ admits an antimagic orientation.
\end{thm}
 \medskip

Following the ideas in the proof of Theorem~\ref{pm} given  in \cite{DTJ}, 
that is, using a different approach to generalizing (a directed version of) Cranston's result \cite{DC}, 
we first prove Theorem~\ref{independentset}  which  generalizes Theorem~\ref{pm}.  

\begin{thm}\label{independentset}
Let $G$ be a graph  and let $A$ be an independent set of $G$ such that   $G$ has a matching  $M$ of size  $ |G|-|A|$ that saturates all vertices in $V(G)\less A$ and $d(v)\geq3$ for every $v\in A$. Then $G$ admits an antimagic orientation.
\end{thm}

We prove   Theorem~\ref{independentset}    in Section~\ref{largealpha}. It is worth noting  that  every  $(a,b)$-biregular bipartite graph with $  a\ge b\ge2 $ has a matching that saturates all vertices in $A$. Theorem~\ref{independentset} implies immediately that every $(a, b)$-biregular bipartite graph with $a\ge b\ge3$ admits an antimagic orientation.   \medskip

 We then study antimagic orientations of graphs $G$ with $\alpha(G)\le4$.  We prove the following Theorem \ref{alpha}.  
 \begin{thm}\label{alpha} Let $G$ be a connected graph.  
\begin{enumerate}[\rm(i)] 
\item If  $\alpha(G)\le2$, then $G$ has an antimagic orientation.
\item If  $\alpha(G)=3$ and $|G|\ge 13$, then $G$ has an antimagic orientation.
\item If  $\alpha(G)=4$ and $\delta(G)\ge11$,  then   $G$ has an antimagic orientation.
\end{enumerate} 
\end{thm}

 Before we prove
 Theorem \ref{alpha}, we first prove a technical result (Theorem \ref{mainresult}) which is instrumental in the proof of Theorem  \ref{alpha}.  
The proof of Theorem~\ref{mainresult} uses the  technique of  Eulerian orientations.     This strategy was previously  used in  \cite{DTJ, LSWY, Y, maxdegree}. 
 Our method here has  new ideas and  is more involved.

 \begin{thm}\label{mainresult}
Let $G$ be a connected graph and let  $X =\{x_1,\ldots,x_t\}\subseteq V(G)$  with $t\in[3]$  such that $d_G(v,X)\le 2 $ for each $v\in V(G)\less X$,  and either 
\begin{enumerate}[\rm(i)] 
\item\label{qtau1}   $t=1$ and $ |N(x_1)|\ge 2 $;  or
 
 \item\label{qtau2}  $t=2$, $e(G)\ge 2|G|-5$   and there exist distinct vertices $y_1, \ldots,y_7\in V(G)\less X $ such that  $x_1$ is complete to either $\{ y_1,y_2, y_3, x_2\}$ or $\{  y_1,y_2, y_3, y_4\}$, and $x_2$ is complete to $\{y_4, \ldots, y_7\}$;  or    
\item\label{qtau3}     $t=3$, $e(G)\ge2|G|-4$ and there exist distinct vertices $y_1, \ldots,y_{11}\in V(G)\less X $ such that  $x_1$ is complete to $\{y_1,y_2, y_3, y_4\}$,  $x_2$ is complete to $\{y_4, y_5, y_6, y_7,   y_8\}$,   and  $x_3$ is complete to either $\{ y_4, y_9,y_{10},y_{11}\}$ or    $\{ y_8, y_9,y_{10},y_{11}\}$.  

 \end{enumerate}
Then $G$ admits an antimagic orientation. 
\end{thm}
 
The proof of Theorem~\ref{alpha} is given in Section~\ref{Mainresult} and the proof of Theorem~\ref{mainresult} is given  in Section~\ref{alpha34}.    For the sake of a cleaner    presentation of the argument, we make no attempt to optimize the constraints in Theorem~\ref{mainresult}.  Corollary~\ref{radius2} follows immediately from Theorem~\ref{mainresult}(i),  and  Corollary~\ref{generalmax}  follows   from Theorem~\ref{dominating} and Corollary~\ref{radius2}. The proof of Corollary~\ref{main} is omitted here as it is similar to  the proof of Theorem~\ref{mainresult}(ii,iii)   when   every vertex of $  V(G)\less X$ is adjacent to a vertex of  $X$ in $G$.

 \begin{cor}\label{radius2}
Every  connected graph $G$ with $r(G)\le 2$  admits an antimagic orientation.

\end{cor}

\begin{cor}\label{generalmax}
Let $G$ be a graph with $|G|\ge 9$. If $\Delta(G)+\delta(G)\ge |G|-2$, then $G$ admits an antimagic orientation.
\end{cor}

 \begin{cor}\label{main}
Let $G$ be a connected graph  and let   $X =\{x_1,\ldots,x_t\}\subseteq V(G)$  with $t\in[3]$  such that  every vertex  in   $    V(G)\less X$ is adjacent to a vertex in $X$,    and either 
\begin{enumerate}[\rm(i)]

 \item\label{qtau1}  $t=2$  and  $\big|N(x_1)\cap N[x_2]\big|\ge1$;  or    
\item\label{qtau2}     $t=3$, $e(G)\ge 2|G|-4$ and there exist distinct vertices  $y_1,y_2\in V(G)\less X$ such that  $y_1$ is complete to $\{x_1,x_2\}$ and $y_2$ is complete to $\{x_2,x_3\}$ in $G$.  

 \end{enumerate}
Then $G$ admits an antimagic orientation. 
\end{cor}


%

     \medskip

\section{Preliminaries}\label{preliminary} 
 In this section,   we shall prove a technical lemma  that  plays an important role in the proofs of our main results.  A closed walk in a multigraph is an \dfn{Euler tour} if
it traverses every edge of the graph exactly once. The following is a result of Euler which shall be needed in the proof of Lemma~\ref{euler}.
\begin{thm}[Euler 1736]\label{Euler}
A connected multigraph admits an  Euler tour if and only if every vertex has even degree.
\end{thm}

 \begin{lem}\label{euler} 
  Let $p\ge0$ be an   integer and let    $G $ be a graph with $m\ge1$ edges. Then there exist an orientation $D $ of $G$ and   bijections $\sigma_i: A(D)\rightarrow   \{p+1,\ldots,p+m\}$,  where $i\in[2]$, such that    for all $v\in V(G)$,  
  \begin{alignat*}{2}
 - \lfloor(d(v)-1)/{2}\rfloor -(p+m) &\le &s_{(D,\sigma_1)}(v) &\le -\lfloor(d(v)-1)/{2}\rfloor+(p+m), \text{ and} \\
 \lfloor(d(v)-1)/{2}\rfloor -(p+m) &\le{} &  s_{(D,\sigma_2)}(v) &\le   \lfloor(d(v)-1)/{2}\rfloor+(p+m).
\end{alignat*}
   \end{lem}
\pf  Let $G$, $m$ and  $p$ be given as  in the statement. We may assume that   $G$ is connected.  Let $A$ be the set (possibly empty) of all vertices $v\in V(G)$ with $d(v)$ odd. Then $|A|=2\ell$ for some integer $\ell\ge0$. Let $G^*:=G$ when $\ell=0$. When $\ell\ge1$, we may assume that $A:=\{x_1, x_2, \ldots, x_{2\ell}\}$.  Let $G^*$ be obtained from $G$ by adding $\ell$ new edges $x_ix_{i+\ell} $ for all $i\in[\ell]$.  Then $e(G^*)=m+\ell$. By Theorem~\ref{Euler},  $G^*$ contains an Euler tour, say $W$,  with vertices and edges $v_1, e_1, v_2, e_2, \ldots v_{m+\ell}, e_{m+\ell}, v_1$ in order, where $v_1, \ldots, v_{m+\ell}$ are not necessarily distinct, and edges $e_1=v_1v_2, e_2=v_2v_3, \ldots, e_{m+\ell}=v_{m+\ell}v_1$ are pairwise distinct. We may further assume that  $e_1\in E(G)$.  Let $1=i_1<i_2<\cdots< i_m\le m+\ell$ be such that $ e_{i_1},\ldots, e_{i_m}$ are all the edges of $G$.  Let $D$ be the orientation of $G$ obtained by orienting  each edge $e_{i_j}$ from $v_{i_j}$ to $v_{{i_j}+1}$ for each  $j\in[m]$, where all arithmetic on indices here and henceforth is done
modulo $m+\ell$.\medskip

Let $\sigma_1: E(G)\rightarrow  \{p+1,\ldots,p+m\}$ be the bijection such that $\sigma_1(e_{i_j})=p+j$ for all $j\in[m]$. It is worth noting  that for all $j, k\in[m]$ with $j<k$,    if $e_{i_j}$ and $e_{i_k}$ are two consecutive edges  on the Euler tour  $W$,  then $i_k=i_j+1$ and  $\sigma_1(e_{i_j})-\sigma_1(e_{i_k})=-1$, that is,   $\sigma_1(e_{i_j})-\sigma_1(e_{i_k})$ contributes $-1$ to   $s_{(D, \sigma_1)}(v_{i_k})$.  By the choice of $G^*$, for each   vertex $v \in V(G)$, there are      at least   $ \lfloor(d(v)-1)/{2}\rfloor $  many such pairs of consecutive edges (incident with $v$)  on the Euler tour $W$.  It follows that  for all $v\in V(G)$,  
\[-\lfloor(d(v)-1)/{2}\rfloor -(p+m)\le s_{(D, \sigma_1)}(v)\le -\lfloor(d(v)-1)/{2}\rfloor +(p+m).\]

Next,  let $\sigma_2: E(G)\rightarrow  \{p+1,\ldots,p+m\}$ be the bijection such that $\sigma_2(e_{i_j})=m+p+1-j$ for all $j\in[m]$.   Then for all $j, k\in[m]$ with $j<k$,    if $e_{i_j}$ and $e_{i_k}$ are two consecutive edges   in the Euler tour  $W$,    then $i_k=i_j+1$ and  $\sigma_2(e_{i_j})-\sigma_2(e_{i_k})=1$, that is,   $\sigma_2(e_{i_j})-\sigma_2(e_{i_k})$ contributes $1$ to   $s_{(D, \sigma_2)}(v_{i_k})$.  It follows that  for all $v\in V(G)$,   
\[\lfloor(d(v)-1)/{2}\rfloor -(p+m)\le s_{(D, \sigma_2)}(v)\le \lfloor(d(v)-1)/{2}\rfloor +(p+m).\]

This completes the proof of Lemma~\ref{euler}.
\qed\bigskip

 \noindent {\bf Remark}.  For the sake of simplicity and clarity of presentation, we  shall  apply  Lemma~\ref{euler} to graphs $H$ with no edges. Under those circumstances, we  shall   let   $D$ with    $V(D)=V(H)$ and   $A(D)=\emptyset$ be the orientation of $H$,   and $\tau: A(D)\rightarrow   \emptyset$  with $s_{(D, \tau)}(v)=0$ for all $v\in V(D)$  be the bijection. \medskip 

  Lemma~\ref{partition} below is a result of Kaplan, Lev and Roditty that  will be needed in the proof of Theorem~\ref{independentset}. 
 
\begin{lem}[\cite{KP}]\label{partition}
Let $t\ge 1$ and $n\ge2$ be   integers and let $n=r_1+\cdots+r_t$ be a partition of $n$, where  $r_i\ge2$ is an integer for all $i\in[t]$.   Then the set $\{1,2,\ldots, n\}$ can be partitioned into  pairwise disjoint subsets  $R_1,   \ldots,R_t$ such that for all $i\in[t]$,  $|R_i|=r_i$,  and $\sum_{r\in R_i}r\equiv0 \emph{ (mod    n+1)}$ if $n$ is even and $\sum_{r\in R_i}r\equiv0 \emph{ (mod    n)}$   if $n$ is odd.  \end{lem}

\section{Proof of Theorem~\ref{independentset}}\label{largealpha}

Let $G$, $A $ and $M$ be given as in the statement.  Let $m:=e(G)$, $n_1:=|A|$ and $n_2:=V(G)\less A$. Then $n_1+n_2=|G|$. By the assumption of $A$ and $M$, we have $ n_1\ge   n_2\ge3$.   Let $A:=\{a_1, \ldots,a_{n_1}\}$ and $V(G)\less A:=\{b_1,\ldots,b_{n_2}\}$ such that $a_ib_i\in M$  for all $i\in[n_2]$.  For each $j\in\{n_2+1, \ldots, n_1\}$,  let  $e_j$ be an arbitrary edge incident with $a_j$ in $G$. Let $E_0:=\{e_{n_2+1},\ldots,e_{n_1}\}$ and let $H:=G\less (M\cup E_0\cup E(G\less A)$. Then $H$ is a bipartite graph with bipartition $(A, V(G)\less A)$.   For all $i\in [n_1]$, since $d_G(a_i)\ge3$, we see that $d_{H}(a_i)\geq2$.     Then $e(H)\ge 2n_1\ge6$  and $e(G\less A)=m-e(H)-|M\cup E_0|=m-n_1-e(H)$. 
  By Lemma \ref{euler}  applied to $G\less A$ with $p=e(H)$, there exist an orientation $D'$ of $G\less A$ and a bijection $\sigma_1:A(D')\rightarrow  \{e(H)+1, \ldots, m-n_1\}$ such that for all  $j\in [n_2]$,   $s_{(D',\sigma_1)}(b_j)\leq m-n_1$.  Let  $D$ be the orientation of $G$ obtained from $D'$ by orienting every edge in $M\cup E_0 \cup E(H)$  towards  $A$.  We next label the edges in $M\cup E_0 \cup E(H)$. \medskip

  For all $i\in[n_1]$, let $A_i$ be the set of all edges incident with $a_i$ in  $H$. Then $|A_i|\ge2$ and $e(H)=|A_1|+\cdots+|A_{n_1}|$. By Lemma~\ref{partition} applied to $e(H)$ with $t=n_1$ and $r_i=|A_i|$ for all $i\in[t]$, the set  $\{1, \ldots, e(H)\}$ can be partitioned into   $R_1, \ldots, R_{n_1}$  such that for all $i\in[n_1]$,  $|R_i|=|A_i|$ and 
$\sum_{r\in R_i}r\equiv0$  (mod   $q$),   where $q= e(H) +1$ if $e(H)$ is even and  $q= e(H)  $ if $e(H)$ is odd. Let $\tau_1: E_0\cup E(H)\rightarrow \{m-|E_0|+1, \ldots, m\}\cup \{1, \ldots, e(H)\}$ be a bijection such that edges in $E_0$ are labelled by integers in $\{m-|E_0|+1, \ldots, m\}$ and for all $i\in[n_1]$,  edges in $A_i$ are labelled by integers in $R_i$.   For each $j\in[n_2]$, let $s_{\tau_1}(b_j)$ denote the sum of labels on all edges incident with $b_j$  in $E_0\cup E(H)$   under $\tau_1$.  We may further assume that 
\[s_{(D',\sigma_1)}(b_1)+s_{\tau_1}(b_1)\ge s_{(D',\sigma_1)}(b_2)+s_{\tau_1}(b_2)\ge\cdots\ge s_{(D',\sigma_1)}(b_{n_2})+s_{\tau_1}(b_{n_2}).\]
 Finally, let $\tau: A(D)\rightarrow [m]$ be the  bijection obtained from $\sigma_1$ and $\tau_1$ by letting $\tau(a_i b_i)=m-n_1+i$ for all $i\in[n_2]$; $\tau(e)=\sigma_1(e)$ for all $e\in E(G\less A)$; and $\tau(e)=\tau_1(e)$ for all $e\in E_0\cup E(H)$. \medskip
 
  By the choice of $(D, \tau)$,  $0>s_{(D,\tau)}(b_1)>  s_{(D,\tau)}(b_2)>\cdots>s_{(D,\tau)}(b_{n_2})$. 
  Since $n_1< e(H)+1$, we see that for all $i, j\in[n_1]$ with $i\ne j$, $m-n_1+i\ne m-n_1+j$ (mod $q$). It follows that  $ s_{(D,\tau)}(a_i)>0$ for all $i\in[n_1]$ and $ s_{(D,\tau)}(a_1), s_{(D,\tau)}(a_2), \ldots, s_{(D,\tau)}(a_{n_1})$ are pairwise distinct. Therefore   $(D, \tau)$ is an antimagic orientation   of $G$.  \medskip

This completes the proof of Theorem~\ref{independentset}. \qed\bigskip

\section{Proof of Theorem~\ref{mainresult}}\label{Mainresult}

Let $G$ and  $X=\{x_1,\ldots,x_t\}$  (and  $y_1, \ldots, y_{4t-1}$ when  $t\in \{2,3\}$)  be given as in the statement. Let   $n:=|G| $,  $m:=e(G)$ and  $\overline{X}:=V(G)\less X$.  Then $m\ge n-1$ and the statement holds for $n\le4$. We may assume that $n\ge5$.  Let
\[A_1:=\big\{v\in \overline{X}: vx_i\in E(G) \text{ for  some } i\in [t]\big\}  \text{, and }  A_2:= \overline{X} \less A_1 \text{ (possibly empty)}.\]
 Then every vertex in $A_1$ is adjacent to at least one vertex in $X$, and
every vertex  $v\in A_2$ is adjacent to at least one vertex in $A_1$ because $d_G(v,X)\le 2 $.  Let $n_i:=|A_i|$ for each $i\in [2]$. Then $n_1+n_2=n-t$.  For each $v\in A_1$, let $e_v:=vx_i$ for some $i\in [t]$ such that when $t=2$, $\{e_{y_1}, \ldots, e_{y_7}\}=\{y_1x_1, y_2x_1, y_3x_1, y_4x_2, y_5x_2, y_6x_2, y_7x_2\}$; and when $t=3$,  $\{e_{y_1}, \ldots,  e_{y_{11}}\}=\{y_1x_1, y_2x_1, y_3x_1, y_4x_2, \ldots, y_8x_2, y_9x_3, y_{10}x_3, y_{11}x_3\}$. Let $E_0:=\{e_v: v\in A_1\}$.
  Let $F$ be a spanning forest of $G$  with $t$ components $F_1,\ldots,F_t$ such that  $E_0\subseteq E(F)$ and for each $i\in[t]$, $x_i\in V(F_i)$   and $d_{F_i}(v,x_i)\leq2$ for every vertex $v\in V(F_i)$.  By the choice of $F$, $N_F(x_1), \ldots, N_F(x_t)$ are pairwise disjoint, $\cup_{i\in[t]}N_F(x_i)=A_1$,  and  for each  $v\in A_2$, $d_F(v)=1$.    Let $H:=G\less E(F)$ and let $E_1$ be the set of all edges  $e$ in $H$ such that  $e$ has an end in  $X$.       Then $E(G[X])\subseteq E_1$ and $m_1:=|E_1|\geq t-1$ for all  $t\in[3]$.  \medskip 
  
 We next show that   there exist  an  orientation $D$ of $G$ and a bijection $\tau:A(D)\rightarrow [m]$  such that   \medskip
 
\myitem{({\bf P1})}  vertices of $A_1$ can be enumerated as $u_1, \ldots, u_{n_1}$ with  $s_{(D,\tau)}(u_1)> \cdots>s_{(D,\tau)}(u_{n_1})>0$; and  
\myitem{({\bf P2})} vertices of $A_2$  can be enumerated as $v_1, \ldots, v_{n_2}$ with  $s_{(D,\tau)}(v_1)< \cdots<s_{(D,\tau)}(v_{n_2})\le0$; and 
\myitem{({\bf P3})}  for each $i\in[t]$,   $s_{(D,\tau)}(x_i)<s_{(D,\tau)}(v_1)$. Furthermore,   when $t=2$,  we have $s_{(D,\tau)}(x_2)+\tau(x_2y_4)<s_{(D,\tau)}(v_1)$;  and  when $t=3$, we have $s_{(D,\tau)}(x_2)+\tau(x_2y_4)+\tau(x_2y_8)<s_{(D,\tau)}(v_1)$. \medskip \\

   \noindent To find such an  orientation $D$, we first orient and label  the edges in $H\less X$. Note that for each $v\in A_1\cup A_2$, $d_{H\less X}(v)\le   n_1+n_2-1$.   By Lemma \ref{euler}  applied to $H\less X$, there exist  an orientation $D'$ of $H\less X$ and   bijections $\sigma_i: A(D') \rightarrow \{m_1+1, \ldots, m-n_1-n_2\}$, where $i\in[2]$,  such that for all $v\in A_1\cup A_2$,
\begin{alignat*}{2}
 - \lfloor(n_1+n_2-2)/{2}\rfloor-(m-n_1-n_2) &\le &s_{(D',\sigma_1)}(v) &\le m-n_1-n_2, \text{ and} \\
-(m-n_1-n_2) &\le{} &  s_{(D',\sigma_2)}(v) &\le    \lfloor(n_1+n_2-2)/{2}\rfloor +m-n_1-n_2.
\end{alignat*}
Let $D$ be the orientation of $G$ obtained from $D'$ by first orienting  all the edges between $A_2$ and $A_1$ in $F$ away from $A_2$,  then all the edges between $X$ and $A_1$ in $G$ away from $X$, and finally    edges of $G[X]$ such that   at most one edge  in $G[X]$ is oriented towards each vertex in $X$.  \medskip

To find such a bijection $\tau$, we next label all the edges in $E_1$ when $E_1\neq \emptyset$.  Let   $\tau_1: E_1   \to [m_1]$  be a bijection  such that   when $t=2$,   $\tau_1(x_1y_4)=1$ if $ y_4\in N(x_1)$ and $\tau_1(x_1x_2)=1$ if $y_4\notin N(x_1)$; and 
 when $t=3$,  $\tau_1(x_1y_4)=1$ and $\tau_1(x_3y_4)=2$ if $y_4\in N(x_3)$,    $\tau_1(x_1y_4)=1$ and $\tau_1(x_3y_8)=2$ if   $y_4\notin N(x_3)$,  and    if $p:=e(G[X])\ge1$, then $\tau_1(E(G[X]))= \{3, \ldots, 2+p\}$. For each $u\in A_1$, let $s_{\tau_1}(u)$ denote the sum of labels on all edges incidents with $u$ in $E_1$ under $\tau_1$.  In the remaining   proof of the existence of $\tau$,  we shall   apply $\sigma_1$ to $A(D')$ when $n_1\le n_2$,  and $\sigma_2$ to $A(D')$ when $n_1>n_2$, that is, we shall apply  $\sigma_1$ and $\sigma_2$  to two separately cases.   Hence when  $A_2\ne\emptyset$, we may  further assume that  $A_2: =\{v_1,v_2,\ldots, v_{n_2}\}$ such that  for each $i\in[2]$, \[s_{(D',\sigma_i)}(v_1)\le s_{(D',\sigma_i)}(v_2)\le \ldots\le s_{(D',\sigma_i)}(v_{n_2}).\] 
For each  $i\in[n_2]$, let $e_i$ be the unique edge incident with $v_i$ in $F$.  Let $E_2:=\{e_i: i\in [n_2]\}$. Then $E_2=E(F)\less E_0$.  We finally label the edges in $E_0\cup E_2$ using integers in $\{m-n_1-n_2+1, \ldots, m\}$ by considering two separate  cases  $n_1\le n_2$ and $n_1> n_2$.\medskip

\medskip
Assume first that  $n_1\le n_2$.  Then  $n_2\ge n_1\ge 2$.  We shall apply $\sigma_1$ to $A(D')$.  To   label the edges in $E_2$, let $\tau_2:    E_2  \to \{m-n_1-n_2+1,\ldots, m- n_1\}$  be the bijection  such that  for each  $i\in[n_2]$,  $\tau_2(e_i)=m-n_1-i+1$.  Let $s_{\tau_2}(u)$ denote the sum of labels on all edges incidents with $u$ in $E_2$ under $\tau_2$ for each $u\in A_1$.
 To label the edges in $E_0$, we need to order the vertices in $A_1$. Let  $A_1: =\{u_1,u_2,\ldots, u_{n_1}\}$ such that   
\[s_{(D',\,\sigma_1)}(u_1)+s_{\tau_1}(u_1)+s_{\tau_2}(u_1)\ge \cdots\ge  s_{(D',\,\sigma_1)}(u_{n_1})+s_{\tau_1}(u_{n_1})+s_{\tau_2}(u_{n_1}).\] 
Finally, let   $\tau: A(D)\rightarrow [m]$ be the  bijection such that 
 $\tau(e_{u_i})=m- i+1$ for each  $i\in [n_1]$, $\tau(e)=\sigma_1(e)$ for each $e\in E(H\less X)$,   $\tau(e)=\tau_1(e)$ for each $e\in E_1$,   and $\tau(e)=\tau_2(e)$ for each $e\in E_2$.  
 Then   for all  $i\in[n_2]$,
\[ s_{(D,\tau)}(v_i)=s_{(D',\, \sigma_1)}(v_i)-\tau_2(e_i) \le (m-n_1-n_2)-(m-n_1-i+1)=-(n_2-i+1)<0,\]
and by the choice of $(D, \tau)$, 
\[
\begin{split}
0>s_{(D,\tau)}(v_{n_2})>\cdots>s_{(D,\tau)}(v_1)&=s_{(D',\, \sigma_1)}(v_1)-\tau_2(e_1)\\
&\ge -\left\lfloor\frac{n_1+n_2-2}{2}\right\rfloor-(m-n_1-n_2)-(m-n_1)\\
&\ge -2(m-n_1)+1-\left\lfloor\frac{n_1-n_2}{2}\right\rfloor \\
&\ge -2(m-n_1)+1\\
 \end{split}
 \]
because $2\le n_1\le n_2$. This proves that $(D, \tau)$ satisfies ({\bf P2}). 
Next for each $j\in[n_1]$, by the orientation of $D$ and the choice of $\tau_1, \tau_2$, we see that $s_{\tau_1}(u_j)\ge0$,  $s_{\tau_2}(u_j)\ge0$, and 
\[\begin{split}
 s_{(D,\tau)}(u_j)&=s_{(D',\,\sigma_1)}(u_j)+s_{\tau_1}(u_j)+s_{\tau_2}(u_j)+\tau(e_{u_j}) \\
 &\ge s_{(D',\,\sigma_1)}(u_j)+\tau(e_{u_j})\\
 &\ge -\left\lfloor\frac{n_1+n_2-2}{2}\right\rfloor-(m-n_1-n_2)+(m-n_1+1) \\
 &\ge  -\left\lfloor\frac{n_1-n_2}{2}\right\rfloor+2\\
 &\geq2, 
  \end{split}
 \]
because $n_1\le n_2$. Hence $s_{(D,\tau)}(u_1)>s_{(D,\tau)}(u_2)>\cdots>s_{(D,\tau)}(u_{n_1})\geq2$. This proves that $(D, \tau)$ satisfies ({\bf P1}).  To see  $(D, \tau)$ satisfies ({\bf P3}), for $x_i\in X$ and each $v\in N_F(x_i)$, $\tau(x_iv)\in   \{m-n_1+1, \ldots,m\}$ because $x_iv\in E_0$. When  $t=1$, then  $d(x_1)\ge 2$ and so $s_{(D,\tau)}(x_1)\le -(m+m-1)=-2m+1 <  s_{(D,\tau)}(v_1)$. When $t\geq2$, then for each $x_i\in X$,   $ d_F(x_i)\ge 3$  and  at most one edge in $G[X]$ is oriented towards $x_i$ in $D$, thus 
\[s_{(D,\tau)}(x_i)\le -(m-n_1+1)-(m-n_1+2)-(m-n_1+3)+4=-3(m-n_1)-2 <s_{(D,\tau)}(v_1).\]
 Finally, when $t=2$,    we have 
 \[
\begin{split}
s_{(D,\tau)}(x_2)+\tau(x_2y_4)&\le 4- \tau(x_2y_5)- \tau(x_2y_6)- \tau(x_2y_7) \\
 &\le   -(m-n_1+1) -(m-n_1+2)-(m-n_1+3)+4   \\
 &<s_{(D,\tau)}(v_{1}).
 \end{split}
 \]
  When $t=3$,   
 we have  
 \[
\begin{split}
s_{(D,\tau)}(x_2)+\tau(x_2y_4) +\tau(x_2y_8)&\le  4- \tau(x_2y_5)- \tau(x_2y_6)- \tau(x_2y_7)  \\
&\le   -(m-n_1+1) -(m-n_1+2)-(m-n_1+3)+4   \\
 &<s_{(D,\tau)}(v_{1}).
 \end{split}
 \]  This proves that $(D, \tau)$ satisfies ({\bf P3}). \medskip

Assume next that  $n_1> n_2$. We may assume that $n_2>0$ because ({\bf P2})  is  trivially true when $n_2=0$.   We shall apply $\sigma_2$ to $A(D')$.   To   label the edges in $E_2$, let $\tau_2:    E_2  \to \{m- n_2+1,\ldots, m \}$  be the bijection  such that  for each  $i\in[n_2]$,  $\tau_2(e_i)=m-i+1$.  Let $s_{\tau_2}(u)$ denote the sum of labels on all edges incidents with $u$ in $E_2$ under $\tau_2$ for $u\in A_1$.
  To label the edges in $E_0$, we need to order the vertices in $A_1$. Let  $A_1: =\{u_1,u_2,\ldots, u_{n_1}\}$ such that   
\[s_{(D',\,\sigma_2)}(u_1)+s_{\tau_1}(u_1)+s_{\tau_2}(u_1)\ge \cdots\ge  s_{(D',\,\sigma_2)}(u_{n_1})+s_{\tau_1}(u_{n_1})+s_{\tau_2}(u_{n_1}).\] 
Finally, let   $\tau: A(D)\rightarrow [m]$ be the  bijection such that 
 $\tau(e_{u_i})=m- n_2-i+1$ for each  $i\in [n_1]$, $\tau(e)=\sigma_2(e)$ for each $e\in E(H\less X)$,   $\tau(e)=\tau_1(e)$ for each $e\in E_1$,   and $\tau(e)=\tau_2(e)$ for each $e\in E_2$.  
 Then   for all  $i\in[n_2]$,
\[ s_{(D,\tau)}(v_i)=s_{(D',\,\sigma_2)}(v_i)-\tau_2(e_i) \le  \left\lfloor\frac{n_1+n_2-2}{2}\right\rfloor+(m-n_1-n_2)-(m- n_2+1)= -2 +\left\lfloor\frac{n_2-n_1}{2}\right\rfloor<0,\]
because $n_1>n_2$. By the choice of $(D, \tau)$, 
\[
\begin{split}
0>s_{(D,\tau)}(v_{n_2})>\cdots>s_{(D,\tau)}(v_1)&=s_{(D',\,\sigma_2)}(v_1)-\tau_2(e_1)\\
&\ge -(m-n_1-n_2)-m\\
&\ge -2 m +n_1+n_2.  \\
 \end{split}
 \]
  This proves that $(D, \tau)$ satisfies ({\bf P2}). 
Next for each $j\in[n_1]$, by the orientation of $D$ and the choice of $\tau_1, \tau_2$, we see that $s_{\tau_1}(u_j)\ge0$,  $s_{\tau_2}(u_j)\ge0$, and 
\[\begin{split}
 s_{(D,\tau)}(u_j)&=s_{(D',\,\sigma_2)}(u_j)+s_{\tau_1}(u_j)+s_{\tau_2}(u_j)+\tau(e_{u_j}) \\
 &\ge s_{(D',\,\sigma_2)}(u_j)+\tau(e_{u_j})\\
 &\ge  -(m-n_1-n_2)+(m-n_2-j+1) \\
 &\ge  n_1-j+1 \\
 &\geq1. 
  \end{split}
 \]
  Hence $s_{(D,\tau)}(u_1)>s_{(D,\tau)}(u_2)>\cdots>s_{(D,\tau)}(u_{n_1})\geq1$. This proves that $(D, \tau)$ satisfies ({\bf P1}).  When  $t=1$, then  $d(x_1)=n_1\ge 3$ and so $s_{(D,\tau)}(x_1)< - (m-n_2)-(m-n_2-1) 
  \le-2m+n_1+n_2   \le s_{(D,\tau)}(v_1)$, so $(D, \tau)$ satisfies ({\bf P3}). To prove  $(D, \tau)$ satisfies ({\bf P3}) when $t\geq2$, since $m\ge 2n-5\ge 2n_1+2n_2-1$, we see that 
  \[4-\sum_{j=1}^3(m-n_1-n_2+j)   = (-2m+n_1+n_2) - (m-2n_1-2n_2+2) < -2m+n_1+n_2  \le s_{(D,\tau)}(v_1).\]
   Then for each $x_i\in X$,   $ d_F(x_i)\ge 3$  and  at most one edge in $G[X]$ is oriented towards $x_i$ in $D$, thus  $s_{(D,\tau)}(x_i)\le -(m-n_1-n_2+1)-(m-n_1-n_2+2)-(m-n_1-n_2+3)+4 <s_{(D,\tau)}(v_1)$.    Finally, when $t=2$,  we have 
 \[
\begin{split}
s_{(D,\tau)}(x_2)+\tau(x_2y_4)&\le 4- \tau(x_2y_5)- \tau(x_2y_6)- \tau(x_2y_7)   \\
&\le  -(m-n_1-n_2+1)  -(m-n_1-n_2+2) -(m-n_1-n_2+3)+4   \\
&<s_{(D,\tau)}(v_{1}).
 \end{split}
 \]
  When $t=3$,   we have 
 \[
\begin{split}
s_{(D,\tau)}(x_2)+\tau(x_2y_4) +\tau(x_2y_8)&\le 4- \tau(x_2y_5)- \tau(x_2y_6)- \tau(x_2y_7)  \\
&\le  -(m-n_1-n_2+1)  -(m-n_1-n_2+2) -(m-n_1-n_2+3)+4   \\
&<s_{(D,\tau)}(v_{1}).
 \end{split}
 \]
Hence $(D, \tau)$ satisfies ({\bf P3}).  This completes the proof  that there exist  an  orientation $D$ of $G$ and a bijection $\tau:A(D)\rightarrow [m]$ such that $(D, \tau)$ satisfies ({\bf P1}), ({\bf P2}) and ({\bf P3}).  \medskip

Since $(D, \tau)$ satisfies  ({\bf P1}), ({\bf P2}) and ({\bf P3}), we see that  $(D, \tau)$ is an antimagic orientation of $G$ if $s_{(D,\tau)}(x_1), \ldots,  s_{(D,\tau)}(x_t)$ are pairwise distinct. We may assume that $s_{(D,\tau)}(x_i)=s_{(D,\tau)}(x_j)$ for some $i, j\in [t]$ with $i\ne j$. Then $t\ge2$ and $E_1\neq \emptyset$. Let $s_i:=s_{(D, \tau)}(x_i)$ for  each $i\in[t]$. Next we will find an antimagic orientation   of $G$ from $(D, \tau)$ by either reversing the direction of an edge in $D[X]$  or strategically swapping the labels on some edges in $E_0\cup E_1$. \medskip

We first  consider the case that $y_4$ is complete to $X$ in $G$. By the choice of $\tau_1$ and $\tau$, we see that $\tau(x_1y_4)=1$,  and $\tau(x_3y_4)=2$ when $t=3$. Let $a:=\tau(x_2y_4)$. 
For $i, j\in[t]$ with $i<j$, let $\tau_{i, j}:A(D)\to [m]$ be the bijection obtained  from $\tau$ by letting $\tau_{i, j}(x_iy_4)=\tau(x_jx_4)$, $\tau_{i, j}(x_jy_4)=\tau(x_ix_4)$ and $\tau_{i, j}(e)=\tau(e)$ for all $e\in E(G)\less\{x_iy_4,x_jy_4\}$. Then for all $v\in V(G)\less \{x_i, x_j\}$, $s_{(D,\tau_{i,j})}(v)=s_{(D,\tau)}(v)$, and so  $(D,\tau_{i,j})$ satisfies ({\bf P1}) and ({\bf P2}).   By ({\bf P3}),  $ s_2+a  <s_{(D,\tau)}(v_1) =s_{(D,\tau_{i,j})}(v_1) $.      
When $t=2$,  we have $s_1=s_2$.   Then $(D,\tau_{1,2})$ is an antimagic orientation of $G$, because $(D,\tau_{1,2})$ satisfies ({\bf P1}) and ({\bf P2}),   $s_{(D,\tau_{1,2})}(x_2) =s_2+a-1 < s_{(D,\tau_{1,2})}(v_{1})$, and     
   \[
\begin{split}
s_{(D,\tau_{1,2})}(x_1)=s_1-(a-1)<s_1+a-1=s_2+a-1=s_{(D,\tau_{1,2})}(x_2)<s_2+a <   s_{(D, \tau_{1,2})}(v_{1}). 
 \end{split}
 \] 
When $t=3$, since $x_1y_4, x_3y_4\in E_1$, so $a\ge 3$.
If $s_1=s_2\ge s_3$,  or $s_1=s_3<s_2$, or $s_2=s_3<s_1$ and $s_1\neq s_2+a-2$, then  $(D,\tau_{2,3})$ is an antimagic orientation of $G$, because  $(D,\tau_{2,3})$ satisfies ({\bf P1}) and ({\bf P2}),      $s_{(D,\tau_{2,3})}(x_2)=s_2+a-2 < s_{(D,\tau_{2,3})}(v_1)$,  and  
  \[
\begin{split}
s_{(D, \tau_{2,3})}(x_3)=s_3-(a-2)<s_1=s_{(D,\tau_{2,3})}(x_1)\ne s_2+a-2= s_{(D,\tau_{2,3})}(x_2) < s_{(D,\tau_{2,3})}(v_1).   
 \end{split}
 \]
Next, if  $s_1= s_2<s_3$,  or $s_1= s_3>s_2$ and $s_2\neq s_1-1$, then $(D,\tau_{1,3})$ is an antimagic orientation of $G$, because  $(D,\tau_{1,3})$ satisfies ({\bf P1}) and ({\bf P2}),   
$s_{(D, \tau_{1,3})}(x_1)= s_1-1< s_{(D,\tau_{1,3})}(v_{1})$, $s_{(D, \tau_{1,3})}(x_2)= s_2$, $s_{(D,\tau_{1,3})}(x_3)= s_3+1$, and so   $s_{(D, \tau_{1,3})}(x_1), s_{(D, \tau_{1,3})}(x_2), s_{(D, \tau_{1,3})}(x_3)$ are pairwise distinct.  Finally, if $s_2=s_3>s_1$, or $s_2=s_3<s_1$ and $s_1= s_2+a-2$, or
 $s_1= s_3>s_2$ and $s_2= s_1-1$,  then  $(D,\tau_{1,2})$ is an antimagic orientation of $G$, because  $(D,\tau_{1,2})$ satisfies ({\bf P1}) and ({\bf P2}), $a\ge 3$,          $s_{(D, \tau_{1,2})}(x_1)= s_1 -(a-1)$, $s_{(D, \tau_{1,2})}(x_2)= s_2+(a-1)< s_{(D,\tau_{1,2})}(v_1)$,   $s_{(D, \tau_{1,2})}(x_3)= s_3$, and so   $s_{(D, \tau_{1,2})}(x_1), s_{(D, \tau_{1,2})}(x_2), s_{(D, \tau_{1,2})}(x_3)$ are pairwise distinct. This completes the proof of the case that $y_4$ is complete to $X$ in $G$.  \medskip

We next   consider the case that $y_4$ is not complete to $X$ in $G$. 
Assume first that   $t=2$. Then    $s_1= s_2$ and  $y_4\notin N(x_1)$. Thus     $x_1x_2\in E(G)$, and  $\tau(x_1x_2)=1$ by the choice of $\tau_1$ and $\tau$. We may assume that the edge $x_1x_2$ is oriented away from $x_1$ in $D$.  Let $D^*$ be obtained from $D$ by    reorienting  the edge $x_1x_2$  towards $x_1$. Then  for all $v\in A_1\cup A_2$, $s_{(D^*,\tau)}(v)=s_{(D,\tau)}(v)$,   and 
\[s_{(D^*,\tau)}(x_2)=s_2-2<s_1+2=s_{(D^*,\tau)}(x_1) =s_2+2  <s_{(D,\tau)}(v_1).\]
  Hence  $(D^*,\tau)$ is an antimagic orientation of $G$. 
  
  Assume next that   $t=3$. Then  $y_4\notin N(x_3)$ and $x_3y_8\in E(G)$.  By the choice of $\tau_1$ and $\tau$,  we have $\tau(x_1y_4)=1$   and $\tau(x_3y_8)=2$.  
Let $a:=\tau(x_2y_4)$ and $b:=\tau(x_2y_8)$. Then $a,b\ge 3$ and $a,b\in \{m-n+4,\ldots, m\}$ are distinct because $x_2y_4, x_2y_8\in E_0$. Note that   $b-2a\le m-2(m-n +4)=-m+2n -8<   -3$ and $a-2b< -3$.   
For each $i \in[2]$, let $\tau_{i, i+1}:A(D)\to [m]$ be the bijection obtained  from $\tau$ by letting $\tau_{i, i+1}(x_iy_{4i})=\tau(x_{i+1}x_{4i})$, $\tau_{i, i+1}(x_{i+1}y_{4i})=\tau(x_ix_{4i})$ and $\tau_{i, i+1}(e)=\tau(e)$ for all $e\in E(G)\less\{x_iy_{4i},x_{i+1}y_{4i}\}$.  Then for all $v\in V(G)\less \{x_i, x_{i+1}\}$, $s_{(D,\tau_{i,i+1})}(v)=s_{(D,\tau)}(v)$, and so  $(D,\tau_{i,i+1})$ satisfies ({\bf P1}) and ({\bf P2}).   By ({\bf P3}),  $ s_2+a +b <s_{(D,\tau)}(v_1)=s_{(D,\tau_{i, i+1})}(v_1) $.   If  $s_2=s_3\ge s_1$,   or  $s_1=s_2< s_3$ and  $s_3\ne s_2+a-1$, or   $s_1=s_3>s_2$ and $s_2+a-1$ is  neither $s_3$  nor  $s_1-(a-1)$, then $(D, \tau_{1,2})$ is an antimagic orientation of $G$, because  $(D,\tau_{1,2})$ satisfies ({\bf P1}) and ({\bf P2}),          $s_{(D, \tau_{1,2})}(x_1)= s_1 -(a-1)$, $s_{(D, \tau_{1,2})}(x_2)= s_2+(a-1)<s_2+a+b<s_{(D,\tau_{1,2})}(v_1)$,   $s_{(D, \tau_{1,2})}(x_3)= s_3$, and so   $s_{(D, \tau_{1,2})}(x_1), s_{(D, \tau_{1,2})}(x_2), s_{(D, \tau_{1,2})}(x_3)$ are pairwise distinct.
Next, if  $s_1=s_2>s_3$, or $s_1=s_3<s_2$, or $s_2=s_3<s_1$ and $s_1\ne s_2+b-2$, 
then $(D, \tau_{2,3})$ is an antimagic orientation of $G$, because  $(D,\tau_{2,3})$ satisfies ({\bf P1}) and ({\bf P2}),          $s_{(D, \tau_{2,3})}(x_1)= s_1$, $s_{(D, \tau_{2,3})}(x_2)= s_2+(b-2)<s_2+a+b<s_{(D,\tau)}(v_{1})=s_{(D,\tau_{2,3})}(v_1)$,   $s_{(D, \tau_{2,3})}(x_3)= s_3-(b-2)$, and so   $s_{(D, \tau_{2,3})}(x_1), s_{(D, \tau_{2,3})}(x_2), s_{(D, \tau_{2,3})}(x_3)$ are pairwise distinct.  \medskip

We now consider the cases that 
   $s_1=s_2 < s_3$ and $s_3=s_2+ a-1$, or   $s_2=s_3<s_1$ and $s_1= s_2+b-2$, or   $s_1=s_3>s_2$, $s_2+a-1 \in\{s_3, s_1-(a-1)\}$ and $b\neq a+1$. 
   Let $\tau^*:A(D)\to [m]$ be the bijection obtained from $\tau$  by  letting $\tau^*(x_1y_4)=a$,   $\tau^*(x_2y_4)=1$, $\tau^*(x_2y_8)=2$,    $\tau^*(x_3y_8)=b$ and  $\tau^*(e)= \tau(e)$ for all $e\in E(G)\less \{x_1y_4, x_2y_4, x_2y_8, x_3y_8\}$. Then for all $v\in V(G)\less X$, $s_{(D,\tau^*)}(v)=s_{(D,\tau)}(v)$, and so  $(D,\tau^*)$ satisfies ({\bf P1}) and ({\bf P2}).   By ({\bf P3}),  $ s_2+a +b <s_{(D,\tau)}(v_1)=s_{(D,\tau^*)}(v_1) $.   Note that  $s_{(D, \tau^*)}(x_1)= s_1-(a-1)$, $s_{(D, \tau^*)}(x_2)=s_2+(a-1)+(b-2)< s_{(D,\tau^*)}(v_1)$,  $s_{(D,\tau^*)}(x_3)= s_3-(b-2)$, and so 
   $s_{(D, \tau^*)}(x_1), s_{(D, \tau^*)}(x_2), s_{(D, \tau^*)}(x_3)$ are pairwise distinct by the assumption and the fact that $a\ge 3$, $b\ge 3$, $b-2a< -3$ and $a-2b<-3$.   \medskip 

It remains to consider the case that   $s_1=s_3>s_2$,  $s_2+a-1 \in\{s_1, s_1-(a-1)\}$ and $b= a+1$. By the choice of $\tau$ and ({\bf P1}),  
\[s_{(D, \tau)}(u_1)- \tau(e_{u_1})\ge s_{(D, \tau)}(u_2) -\tau(e_{u_2})\ge\cdots\ge s_{(D, \tau)}(u_{n_1}) -\tau(e_{u_{n_1}}).\]  Since $b=a+1$, $a=\tau(e_{y_4})$ and $b=\tau(e_{y_8})$,  we see that    $s_{(D, \tau)}(y_4)- a\le s_{(D, \tau)}(y_8)- b$, that is, $s_{(D, \tau)}(y_4)\le s_{(D, \tau)}(y_8)- 1$. We may assume that $y_8=u_{\ell}$ and $y_4=u_{\ell+1}$ for some $\ell\in [n_1-1]$.  We may further assume that $s_{(D, \tau)}(y_4)<  s_{(D, \tau)}(y_8)- 1$  (else,    let   $u_\ell:=y_4$,  $u_{\ell+1}:=y_8$, $\tau(x_2y_4)=b$ and $\tau(x_2y_8)=a$, we have  $s_1=s_3>s_2$  but  $\tau(x_2y_8)\ne \tau(x_2y_4)+1$. By the previous  cases, we see that $G$ admits an antimagic orientation).   
Assume that $s_{(D, \tau)}(y_4)\le  s_{(D, \tau)}(y_8)-3$. Let $\tau^*:A(D)\to [m]$ be the bijection obtained from $\tau$  by  letting $\tau^*(x_1y_4)=2$,   $\tau^*(x_3y_8)=1$,  and  $\tau^*(e)= \tau(e)$ for all $e\in E(G)\less \{x_1y_4, x_3y_8\}$. 
Then for all $v\in V(G)\less \{x_1, x_3, y_4, y_8\}$, $s_{(D,\tau^*)}(v)=s_{(D,\tau)}(v)$, and so  $(D,\tau^*)$ satisfies  ({\bf P2}).    Note that   $s_{(D, \tau^*)}(y_4)=s_{(D,\tau)}(y_4)+1 \le  (s_{(D,\tau)}(y_8)-3)+1< s_{(D,\tau)}(y_8)-1=s_{(D,\tau^*)}(y_8)$. 
It follows that \[s_{(D,\tau^*)}(u_1)>\cdots>s_{(D,\tau^*)}(u_{\ell-1})>s_{(D,\tau^*)}(y_8)>s_{(D,\tau^*)}(y_4)>s_{(D,\tau^*)}(u_{\ell+2})>\cdots>s_{(D,\tau^*)}(u_{n_1}),\]
    and $s_{(D,\tau^*)}(x_1)=s_1-1$, $s_{(D,\tau^*)}(x_2)=s_2$, $s_{(D,\tau^*)}(x_3)=s_1+1< s_2+a+b<s_{(D,\tau^*)}(v_1)$ are pairwise distinct because $s_2+a-1 \in\{s_1, s_1-(a-1)\}$. Hence $(D, \tau^*)$ is an antimagic orientation of $G$.  
Finally, we may assume that  $s_{(D, \tau)}(y_4)=  s_{(D, \tau)}(y_8)-2$.    Let $\tau^*:A(D)\to [m]$ be the bijection obtained from $\tau$  by    letting $\tau^*(x_1y_4)=2$,  $\tau^*(x_3y_8)=1$, $\tau^*(x_2y_4)=b=a+1$,  $\tau^*(x_2y_8)=a$,  and  $\tau^*(e)= \tau(e)$ for all $e\in E(G)\less \{x_1y_4, x_2y_4, x_2y_8, x_3y_8\}$. 
Then for all $v\in V(G)\less \{x_1,x_2, x_3, y_4, y_8\}$, $s_{(D,\tau^*)}(v)=s_{(D,\tau)}(v)$, and so  $(D,\tau^*)$ satisfies  ({\bf P2}).    Note that   $s_{(D, \tau^*)}(y_4)=s_{(D,\tau)}(y_4)+(b+2)-(a+1)=  (s_{(D,\tau)}(y_8)-2)+2= s_{(D,\tau)}(y_8)$, and $s_{(D,\tau^*)}(y_8)=s_{(D,\tau)}(y_8)+(a+1)-(b+2)=s_{(D,\tau)}(y_8)-2=s_{(D,\tau)}(y_4)$. 
It follows that \[s_{(D,\tau^*)}(u_1)>\cdots>s_{(D,\tau^*)}(u_{\ell-1})>s_{(D,\tau^*)}(y_4)>s_{(D,\tau^*)}(y_8)>s_{(D,\tau^*)}(u_{\ell+2})>\cdots>s_{(D,\tau^*)}(u_{n_1}),\]
    and $s_{(D,\tau^*)}(x_1)=s_1-1$, $s_{(D,\tau^*)}(x_2)=s_2$, $s_{(D,\tau^*)}(x_3)=s_1+1< s_2+a+b<s_{(D,\tau^*)}(v_1)$ are pairwise distinct because $s_2+a-1 \in\{s_1, s_1-(a-1)\}$. Hence $(D, \tau^*)$ is an antimagic orientation of $G$.  \medskip
 
This completes the proof of Theorem \ref{mainresult}.\qed\bigskip

\section{Proof of Theorem~\ref{alpha}}\label{alpha34}
   Let $G$  be a connected graph with $\alpha:=\alpha(G)\le4$. By Theorem~\ref{dominating} and the fact that $G$  is connected, Theorem \ref{alpha}(i) holds. To prove Theorem \ref{alpha}(ii,iii), let  $n:=|G|$,  $m:=e(G)$,  and     $S:=\{u_1,  \ldots,  u_\alpha\}$ be a maximum independent set in $G$ such that $H_S$ has the minimum number of components, where $H_S$ is the bipartite subgraph of $G$ with bipartition $\{S, \overline S\}$ and $E(H_S)$ consisting of all edges between $S$ and  $\overline S:=V(G)\less S$ in $G$.  Then every vertex in $\overline S$ is adjacent to at least one vertex in $S$, and each component of $H_S$ contains at least one vertex of $S$, because $\alpha(G)=\alpha$.   By Corollary~\ref{radius2},  we may   assume that  $r(G)\ge3$. Then no vertex in $\overline S$ is complete to $S$ in $G$.  
  We next show that\medskip 

\noindent ($\ast$) every component of $H_S$ contains at least two vertices in $S$. \medskip

\pf Suppose for a contradiction that    $H_1$ is a component of $H_S$ with  $|V(H_1)\cap S|=1$, say $V(H_1)\cap S=\{u_{\alpha}\}$.  Then $H_S$ is disconnected. Since $G$ is connected, there must exist $x\in N(u_\alpha)$ and $x’ \in \bigcup_{i=1}^{\alpha-1} N(u_i)$ such that $xx’\in E(G)$. Note that $x$ is anti-complete to $\{ u_1,\ldots,u_{\alpha-1}\}$ in $G$.  But then $S’:=\{x,u_1,\ldots,u_{\alpha-1}\}$ is a maximum independent set of $G$ such that $H_{S’}$ has fewer components than $H_S$, contrary to the choice of $S$. \qed\medskip

  To prove Theorem~\ref{alpha}(ii), assume $\alpha=3$ and  $|G|\ge 13$. Then  $\overline{G}$ is $K_4$-free. By Tur\'an's Theorem, $e(\overline{G})\le n^2/3$. Hence $m \ge {n\choose 2}-n^2/3   >2n-5$ because $n\ge 13$. Since every component of $H_S$ contains at least two vertices in $S$ and $r(G)\ge3$, we see that $H_S$ is connected. This, together with the fact that $H_S$ is bipartite,   implies that    there must exist distinct vertices $  y_1, y_2 \in \overline S$   such that   $y_1$ is complete to, say  $\{u_1, u_2\}$, and $y_2$ is complete to, say  $\{u_2, u_3\}$ in $H_S$. Then $y_1\ne y_2$ because   $r(G)\ge3$.  By Corollary~\ref{main}(ii)  applied to $G$ with $X=S$, $G$ admits an antimagic orientation.\medskip

  To prove Theorem~\ref{alpha}(iii), assume  $\alpha =4$ and $\delta(G)\ge11$.  Then $n\ge 12$ and $m >5n$. Since every component of $H_S$ contains at least two vertices in $S$, we see that $H_S$ has at most two components. We first consider the case that $H_S$ is connected. Assume first there exists a vertex $x_1\in \overline S$ such that $x_1$ is adjacent to  three vertices of $S$, say $u_1,u_2,u_3$. Then  $x_1u_4\notin E(G)$ because $r(G)\ge3$. Let $x_2\in N(u_4)\cap \big(N(u_1)\cup N(u_2)\cup N(u_3)\big)$, say $x_2\in N(u_4)\cap N(u_3)$. This is possible because $H_S$ is connected. Then $S\subseteq N(x_1)\cup N( x_2)$  and for all $v\in V(G)\less\{x_1, x_2\}$, we have $d(v, \{x_1, x_2\})\le2$.  Note that $\delta(G)\ge11$. 
  By Theorem \ref{mainresult}(\ref{qtau2}) applied to $G$ with $X:=\{x_1,x_2\}$ and $y_4=u_3$,    $G$ admits an antimagic orientation. Next, assume that  no vertex in $\overline S$ is adjacent to three vertices in $S$. Then there exist $x_1,x_2,x_3\in \overline S$ such that  for each $i\in [3]$, $x_i$ is adjacent to, say $\{u_i,u_{i+1}\}$ (and so anti-complete to $S\less \{u_i,u_{i+1}\}$) in $G$. Then $S\subseteq N(x_1)\cup N( x_2)\cup N(x_3)$ and for all $v\in V(G)\less\{x_1, x_2\}$, we have $d(v, \{x_1, x_2, x_3\})\le2$.   
   Since $\delta(G)\ge11$, 
 by Theorem \ref{mainresult}(\ref{qtau3}) applied to $G$ with $X=\{x_1,x_2,x_3\}$, $y_4=u_2$ and  $y_8=u_3$,   $G$ admits an  antimagic orientation.\medskip
  
It remains to consider the case that $H_S$ is disconnected. By ($\ast$),   $H_S$ has exactly two components.  We may further assume that  for every maximum independent set  $S'$ of  $G$, $H_{S'}$ is disconnected and  if $H_{S'}$ has exactly two components,  then   each component of $H_{S'}$ contains exactly two vertices of $S'$ by ($\ast$). Let $H_1, H_2$ be the   components of $H_S$. By ($\ast$),  we may   assume that $u_1,u_2\in V(H_1)$ and $u_3,u_4\in V(H_2)$.   Since $G$ is connected, there must exist $x_1\in N(u_1)\cup N(u_2)$ and $x_2\in N(u_3)\cup N(u_4)$ such that $x_1x_2\in E(G)$. Then $x_1$ is complete to $\{u_1,u_2\}$ and $x_2$ is complete to $\{u_3,u_4\}$ in $G$, else, say  $x_1u_1\notin E(G)$, then $S’:=\{x_1,u_1,u_3,u_4\}$ is a maximum independent set of $G$ such that either $H_{S’}$ is connected,  this contradicts the choice of $S$ or $H_{S’}$ has exactly two components with one component  containing only one vertex $u_1$ in $S’$, contradiction again. Thus $S\subseteq N(x_1)\cup N( x_2)$  and for all $v\in V(G)\less\{x_1, x_2\}$, we have $d(v, \{x_1, x_2\})\le2$.
By Theorem \ref{mainresult}(\ref{qtau2}) applied to $G$ with $X=\{x_1,x_2\}$ and $x_1x_2\in E(G)$, $G$ admits an antimagic orientation. \medskip

    This completes the proof of Theorem~\ref{alpha}.  \qed 

\end{document}